\documentclass[12pt]{amsart}
\usepackage{amsmath}
\usepackage{amssymb}
\usepackage{amsfonts}

\usepackage{hyperref,url}
\usepackage{xcolor}
\hypersetup{
   colorlinks,
    linkcolor={black!50!black},
    citecolor={black!50!black},
    urlcolor={black!80!black}
}

\newcommand{\Z}{\mathbb{Z}}
\newcommand{\N}{\mathbb{N}}

\newcommand{\rseven}{$(r,s)$-even\ }

\newtheorem{thm}{Theorem}[section]
\newtheorem{cor}[thm]{Corollary}

\newtheorem{rem}[thm]{Remark}

\title[Restricted linear congruence]{On the number of solutions of a restricted linear congruence}
\author[K V Namboothiri]{K Vishnu Namboothiri}
\address{Department of Mathematics, Government Polytechnic College, Vennikkulam, Thiruvalla, Kerala - 689 544, INDIA\\Department of Collegiate Education, Government of Kerala, INDIA}
\email{kvnamboothiri@gmail.com}

\begin{document}
\baselineskip=17pt

\begin{abstract}
Consider the linear congruence equation $${a_1^{s}x_1+\ldots+a_k^{s} x_k \equiv b\,(\text{mod } n^s)}\text { where } a_i,b\in\Z,s\in\N$$
Denote by $(a,b)_s$ the largest $l^s\in\N$ which divides $a$ and $b$ simultaneously. Given $t_i|n$, we seek solutions $\langle x_1,\ldots,x_k\rangle\in\Z^k$ for this linear congruence with the restrictions $(x_i,n^s)_s=t_i^s$.  Bibak \emph{et al.} [J. Number Theory,  171:128-144, 2017] considered the above linear congruence with $s=1$ and gave a formula for the number of solutions in terms of the  Ramanujan sums. In this paper, we derive a formula for the number of solutions of the above congruence for arbitrary $s\in\N$ which involves the generalized Ramanujan sums defined by E. Cohen [Duke Math. J, 16(85-90):2, 1949] 
\end{abstract}



\maketitle

\section{Introduction}
Let $a_i,b\in\Z$ for $i=1,\ldots k$ and $n,s\in\N$. Consider the linear congruence equation 
\begin{equation}\label{eq:gen_lin_cong}
{a_1^{s}x_1+\ldots+a_k^{s} x_k \equiv b\,(\text{mod } n^s)}
\end{equation}
  Such equations were considered by many authors who attempted find either their solutions or the number of solutions. For the case $s=1$, such an attempt was made by D. N. Lehmer \cite{lehmer1913certain} who proved that
\begin{thm}
 Let $a_1,\ldots,a_k,b,n\in\Z, n\geq 1$. The linear congruence equation $${a_1 x_1+\ldots+a_k  x_k \equiv b\,(\text{mod } n )}$$ has a solution $\langle x_1,\ldots, x_n\rangle \in \Z_n^k$ if and only if $l|b$ where $l$ is the gcd of $a_1,\ldots,a_k,n$. Furthermore, if this condition is satisfied, then there are $ln^{k-1}$ solutions.
\end{thm}

Later the focus was shifted to solving the above type of congruence equations with some extra restrictions on the solutions $x_i$. One such restriction is requiring that the solutions should satisfy $gcd(x_i,n)=t_i \,(1\leq i \leq k)$ where $t_i$ are given positive divisors of $n$. A linear congruence with such  restrictions is called to be a restricted linear congruence. These kind of restricted congruences  were tried to solve by many authors some of which were special cases of the problem we are going to solve here. With $a_i=s=1$ and restrictions $(x_i,n)=1$, Rademacher \cite{radekacher1925aufgabe}  and Brauer \cite{brauer1926losung} independenty gave a formula for the number of solutions $N_n(k,b)$ of the congruence. An equivalent formula involving the Ramanjans sums was proved by Nicol and Vandiver \cite{nicol1954sterneck}, and E. Cohen \cite{cohen1955class}. The number of solutions given by them is
\begin{equation}
 N_n(k,b) = \frac{1}{n}\sum\limits_{d|n}c_d(b)\left(c_n\left(\frac{n}{d}\right)\right)^k
\end{equation}
where $c_r(n)$ denote the usual Ramanujan sum.

The restricted congruence (\ref{eq:gen_lin_cong}) (with $s=1$) and their solutions has found interesting applications in various fields including number theory, cryptography, combinatorics, computer science etc.  Liskovets defined a multivariate arithmetic function in \cite{liskovets2010multivariate}. The special case of our restricted congruence problem with $b=0$ and $a_i=s=1$ is related to this multivariate function. This functions has many combinatorial as well as topological applications. In computer science, the restricted congruence problem has applications in studying universal hashing (see Bibak \emph{et al.} \cite{bibak2015almost}).  

In \cite{bibak2016restricted} Bibak \emph{et al.} considered the  linear congruence (\ref{eq:gen_lin_cong}) taking $a_i=s=1$ and the restrictions $(x_i,n)=t_i$ where $t_i$ are given positive divisors of $n$. This was later generalized for an arbitrary $s$ with still requiring $a_i=1$ by K V Namboothiri in \cite{namboothiri4}. It was proved there that

\begin{thm}
 Let $b,n\in\Z, n\geq 1$, and $d_1\ldots, d_{\tau(n)}$ be the positive divisors of $n$. For $1\leq l\leq \tau(n)$, define $\mathcal{C}_{j,s}=\{1\leq x \leq n^s|(x,n^s)_s=d^s_j\}$. The number of solutions of the linear congruence 
 \begin{equation*}\label{res_lin_cong}
x_1+\ldots +x_k\equiv b\,(\text{mod }n^s)
 \end{equation*}
with restrictions $(x_i,n^s)_s=d_i^s$  for $i=1,\ldots,k$ is
\begin{equation}
\frac{1}{n^s}\sum\limits_{d|n}c_{d,s}(b)\prod\limits_{j=1}^{\tau(n)}\left(c_{\frac{n}{d_j},s}(\frac{n^s}{d^s})\right)^{g_j}
\end{equation}
 where $g_j= |\{x_1,\ldots, x_k\}\cap \mathcal{C}_{j,s}|$, $1\leq j\leq \tau(n)$.
\end{thm}

In the above statement, $\tau(n)$ stands for the number of positive divisors of $n$.
This result has been proved in some special cases by many authors, for example in \cite{cohen1955class, dixon1960finite, nicol1954sterneck, sander2013adding}. The above result with $a_i=s=t_i=1$ was proved by Cohen \cite[Theorem 12, Theorem 12']{cohen1956extensionof} and reproved using different techniques by Namboothiri \cite[Theorem 6.5]{namboothiri5}. In \cite{bibak2017restricted}, Bibak \emph{et al.} gave a formula for the number of solutions of the restricted congruence (\ref{eq:gen_lin_cong}) with  $s=1$. Following the methods used there, we give an analogous formula for the  equation (\ref{eq:gen_lin_cong}) where $s$ is  arbitrary. Our formula is is precisely 
\begin{equation*}
\frac{1}{n^s}\prod\limits_{i=1}^k\frac{J_s(\frac{n}{t_i})}{J_s(\frac{n}{t_i d_i})}\sum\limits_{d|n}c_{d,s}(b)\prod\limits_{i=1}^k c_{\frac{n}{t_i d_i},s}(\frac{n}{d}) \text{ where } d_i=(a_i,\frac{n}{t_i})
\end{equation*}
Here $J_s$ is the Jordan totient function which generalizes the Euler totient function $\varphi$.

We employ the tools of discrete Fourier transform of $(r,s)$-even functions and special properties of generalized Ramanujan sums given by Cohen \cite{cohen1949extension} to arrive at our formula. As we have already mentioned above, the formula generalizes the main result in \cite{namboothiri4} and the main theorem in \cite{bibak2017restricted}.

\section{Notations and basic results}
For $a,b\in \Z$ with atleast one of them non zero, the generalized gcd of these numbers $(a,b)_s$ is defined to be the largest $l^s\in\N$ dividing $a$ and $b$ simultaneously. Therefore $(a,b)_1=(a,b)$, the usual gcd of two integers. For a prime $p$ and integer $n$, by $p^r||n$ we mean that $p^r|n$, but $p^{r+1}\nmid n$.
The Jordan totient function $J_s$ is defined by $$J_s(n)=n^s\prod\limits_{p|n}(1-p^{-s})$$
By $e(x)$, we mean the function $e(x)=exp(2\pi ix)$. Let $c_r(n)$ denote the Ramanujan sum which is defined to be the sum of $n^{\text{th}}$ powers of primitive  $r^{\text{th}}$ roots of unity. That is, 
\begin{equation}\label{ram_sum}
 c_r(n)=\sum\limits_{j=1, (j,r)=1}^r e\left(\frac{jn}{r}\right)
\end{equation}
 For a positive integer $s$, E. Cohen \cite{cohen1949extension} generalized the Ramanujan sum defining $c_{r,s}$ as 
\begin{equation}
 c_{r,s}(n)=\sum\limits_{j=1,(j,r^s)_s=1}^{r^s}e\left(\frac{jn}{r^s}\right)
\end{equation}
Note that for $s=1$, this definition gives the usual Ramanujan sum defined in equation (\ref{ram_sum}). In the same paper, Cohen also gave the following formula:
\begin{equation}\label{mob_grs}
 c_{r,s}(n)=\sum\limits_{d|r, d^s|n}\mu\left(\frac{r}{d}\right)d^s
\end{equation}
where $\mu$ is the M\"{o}bius function. A special property of $\mu$ is that
\begin{equation}\label{eq:mobiusproperty}
 \sum\limits_{d|n}\mu(d)=\begin{cases}
                          1 \text{ if }n=1\\
                          0 \text{ otherwise}
                         \end{cases}
\end{equation}

An arithmetic function $f$ is said to be periodic with period $r$ (or $r$-periodic) for $r\in\N$ if for every $m\in\Z$, $f(m+r)=f(m)$. The complex exponential funcion $e(x)$ has period 1.

Let $r,s$ be  positive integers. A function $f$ that satisfies $f(m)=f((m,r^s)_s)$ is called as an $(r,s)$-even function. This concept was introduced by McCarthy in \cite{mccarthy1960generation} and many of its properties were studied there. An $(r,s)$-even function is $r^s$-periodic. $c_{r,s}(n)$ is $(r,s)$-even and so $r^s$-periodic. For a proof of these  statements, see \cite[Section 2]{ namboothiri4}.

For an $r$-periodic arithmetic function $f(n)$, its finite Fourier transform is defined to be the function 
\begin{equation}
 \hat{f}(b) = \frac{1}{r}\sum\limits_{n=1}^r f(n) e\left(\frac{-bn}{r}\right)
\end{equation}
A Fourier representation of $f$ is given by 
\begin{equation}
 f(n)= \sum\limits_{b=1}^r \hat{f}(b) e\left(\frac{bn}{r}\right)
\end{equation}
It further satisfies $\hat{\hat{f}}=rf$.
See, for example, \cite{montgomery2006multiplicative} for a detailed study on Finite Fourier transforms. 

Further, for an $(r,s)$-even function $f$, we have  \begin{eqnarray}
 \hat{f}(n) &=& \sum_{k(\text{mod } r^s)} f(k )e(-\frac{kn}{r^s})\nonumber\\
 &=& \sum_{d|r}f(d^s)\sum\limits_{\substack{1\leq j\leq \frac{r^s}{d^s} \\ (j, \frac{r^s}{d^s})_s = 1}}e(-\frac{jd^sn}{r^s})\nonumber\\
 &=& \sum_{d|r}f(d^s)c_{\frac{r}{d},s}(n)\label{fhat_crs}
\end{eqnarray}
Since $f$ and $c_{\frac{r}{d},s}(.)$ are \rseven, $\hat{f}(n)$ is also \rseven. 

The Cauchy product of two $r$-periodic functions $f$ and $g$ is defined to be $(f\otimes g)(n) = \sum\limits_{\substack{1\leq a,b\leq n\\a+b=n\,(\text{mod }r)}}f(a)g(b)$. The definition can be extended to finitely many functions in the same way Note that the DFT of Cauchy product satisfy $$\widehat{f\otimes g}=\hat{f}\hat{g}$$
See \cite[Chapter 2]{terras1999fourier} for a proof of the above property of the Cauchy product.

We divide the proof of our main theorem into mainly two special cases. One of it proves the case when $k=1$. The second one proves the case $a_i=1$. Note that though the case $a_i=1$ was proved in \cite{namboothiri4}, we here give an alternate formula for the number of solutions.

\section{Number of Solutions of the linear congruence in the single variable case}
Recall that our aim is to derive a formula for the number of solutions of the restricted linear congruence equation (\ref{eq:gen_lin_cong}). We derive a necessary and sufficient condition for solutions to exist in the case $k=1$. We also find the number of solutions in this case. This result generalizes \cite[Theorem 3.1]{bibak2017restricted}. 
\begin{thm}\label{th:singlevar}
 Consider the linear congruence equation 
 \begin{equation}\label{eq:singlevar}
  a^{s} x\equiv b\, (\text{mod }n^s)
   \end{equation}where $a,b\in\Z,n,s\in\N$ are given. Given $t\in\N$ where $t|n$, this congruence equation has a solution $x$ with $(x,n^s)_s=t^s$ if and only if $t^s|(b,n^s)_s$ and $(a^{s},\frac{n^s}{t^s})_s=(\frac{b}{t^s},\frac{n^s}{t^s})_s$. The number of solutions in this case is $$d^s\prod\limits_{\substack{p|d\\p\nmid\frac{n}{dt}\\p\text{ prime}}} (1-p^{-s})= \frac{J_s(\frac{n}{t})}{J_s(\frac{n}{dt})}$$ where $d^s = (a^{s},\frac{n^s}{t^s})_s$.
   \begin{proof}
    Assume that there exists $x_0$ such that $a^{s} x_0\equiv b\, (\text{mod }n^s)$ with $(x_0,n^s)_s=t^s$. Now $(a^{s} x_0, n^s)_s = (b,n^s)_s$. Since $(x_0,n^s)_s=t^s$, $(a^{s} x_0, n^s)_s = d^s t^s$ for some $d\geq 1$. So $(b,n^s)_s = d^st^s$ and so  $t^s|(b,n^s)_s$. Now it follows that $$(a^{s} \frac{x_0}{t^s}, \frac{n^s}{t^s})_s = (\frac{b}{t^s},\frac{n^s}{t^s})_s = d^s$$
    But $(\frac{x_0}{t^s}, \frac{n^s}{t^s})_s = 1$ implies that $(a^{s}, \frac{n^s}{t^s})_s = d^s = (\frac{b}{t^s},\frac{n^s}{t^s})_s $.
    
    Conversely, let  $t^s|(b,n^s)_s$ and $(a^{s},\frac{n^s}{t^s})_s=(\frac{b}{t^s},\frac{n^s}{t^s})_s = d^s$ for some $d\geq 1$. So $(\frac{a^{s}}{d^s},\frac{n^s}{d^s t^s})_s=(\frac{b}{d^s t^s},\frac{n^s}{d^s t^s})_s = 1$. This tells us that $(\frac{a^{s}}{d^s},\frac{n^s}{d^s t^s})=1$, or equivalently, that $\frac{a^{s}}{d^s}$ is invertible modulo $\frac{n^s}{d^s t^s}$. Putting $\frac{n}{d t}=N$, we see that $\frac{a^{s}}{d^s}y = \frac{b}{d^s t^s}(\text{mod }N^s)$ has a unique solution modulo $N^s$. So we get a $y_0$ such that $$\frac{a^{s}}{d^s}y_0 = \frac{b}{d^s t^s}(\text{mod }N^s)$$ where $y_0$ is unique modulo $N^s$. Therefore we have $a^{s}(t^s y_0)\equiv b (\text{mod }n^s)$. So $x_0 = t^s y_0$ is a solution for the linear congruence (\ref{eq:singlevar}).

    Now we will find the number of solutions when the given conditions are satisfied. So we are assuming that $t^s|(b,n^s)_s$ and $(a^{s},\frac{n^s}{t^s})_s=(\frac{b}{t^s},\frac{n^s}{t^s})_s = d^s$ for some $d\geq 1$ and so there exists an $x_0$ such that $a^{s} x_0\equiv b (\text{mod }n^s)$ where $(x_0,n^s)_s=t^s$. Therefore $t^s|x_0$ or $x_0=t^sy_0$ for some $y_0$. It follows that $$\frac{a^{s}}{d^s}y_0 = \frac{b}{d^s t^s}(\text{mod }N^s)$$  where  $N = \frac{n}{dt}$. By assumption, $(\frac{a^{s}}{d^s},N^s)_s=1$. Therefore  $y_0$ above is unique modulo $N^s=\frac{n^s}{d^s t^s}$. So the possibile values of $x_0$ are $t^s(y_0+k N^s) = y_0+k\frac{n^s}{d^s}$ modulo $n^s$ which are unique for $0\leq k\leq d^s -1$. We need to count all such $x_0$ to get the number of solutions of the given congruence equation. In short, we need to count $t^s(y_0+k N^s)$ which are solutions of $a^{s}x\equiv b(\text{mod }N^s)$ with $(t^s(y_0+k N^s), n^s)_s=t^s$. But  
    \begin{eqnarray*}
(t^s(y_0+k N^s), n^s)_s=t^s&\Rightarrow& (y_0+k  N^s, \frac{n^s}{t^s})_s = 1\\
&\Rightarrow& (y_0+k  N^s, d^sN^s)_s = 1\\
&\Rightarrow& (y_0+k  N^s, d^s)_s = 1
    \end{eqnarray*}
    Now $(t^s y_0,n^s)_s=t^s\Rightarrow (y_0,N^s)_s=1$.
    Therefore
$$ (y_0+k  N^s, d^s)_s = 1 \Rightarrow (y_0+k  N^s, N^s d^s)_s = 1$$
So we conclude that $$ (y_0+k  N^s, d^s)_s = 1\Leftrightarrow (y_0+k  N^s, N^s d^s)_s = 1$$
Therefore we need to count $t^s(y_0+k N^s)$  with $(y_0+k N^s, d^s)_s=1$ where $0\leq k\leq N^s$. Let us denote the number of solutions by $S$. Then by the property of the M\"{o}bius function given in (\ref{eq:mobiusproperty}), we have
\begin{eqnarray*}
 S &=&\sum\limits_{0\leq k\leq d^s-1}\sum\limits_{l|(y_0+kN^s,d^s)_s}\mu(l)\\
 &=&\sum\limits_{0\leq k\leq d^s-1}\sum\limits_{\substack{l|\delta^s\\\delta^s|(y_0+kN^s,d^s)_s}}\mu(l)\\
 &=&\sum\limits_{0\leq k\leq d^s-1}\sum\limits_{\substack{l|\delta^s\\\delta^s|y_0+kN^s\\ \delta^s|d^s}}\mu(l)\\
 &=&\sum\limits_{0\leq k\leq d^s-1}\sum\limits_{\substack{\delta^s|y_0+kN^s\\ \delta^s|d^s}}\mu(\delta)\\
 &=&\sum\limits_{\delta|d}\sum\limits_{0\leq k\leq d^s-1}\sum\limits_{\substack{\delta^s|y_0+kN^s\\}}\mu(\delta)\\
 &=&\sum\limits_{\delta|d}\mu(\delta) \sum\limits_{\substack{0\leq k\leq d^s-1\\kN^s = -y_0(\text{mod }\delta^s)}}1\\
\end{eqnarray*}
Now $kN^s = -y_0(\text{mod }\delta^s)$ has solution $k$ if and only if $(N^s,\delta^s)|(-y_0)$. But $(N^s,\delta^s)=v^s=(N^s,\delta^s)_s$ for some $v\geq 1$. We claim that $v=1$. Otherwise, if  $v>1$, then $v^s|(-y_0)\Rightarrow v^s|(y_0,N^s)_s=1$ contradicting our assumption that ($y_0,N^s)_s = 1$. So if $(N^s,\delta^s)>1$, the inner sum in $S$ is empty and so is 0. Hence we need to consider only those $\delta$ for which $(N^s,\delta^s)=1$. But in that case the solution $k$ for $kN^s = -y_0(\text{mod }\delta^s)$ is unique mod $\delta^s$. Then the number of solutions modulo $d^s$ would be $d^s/\delta^s$ (which are $k,k+\delta^s, \ldots, k+(r^s-1)\delta^s$ where $r\delta=d$). To conclude, we have

\begin{eqnarray*}
 S &=&\sum\limits_{\substack{\delta|d\\(N^s,\delta^s)=1}}\mu(\delta) \frac{d^s}{\delta^s}\\
 &=&d^s\sum\limits_{\substack{\delta|d\\(N^s,\delta^s)=1}} \frac{\mu(\delta)}{\delta^s}\\
 &=& d^s\prod\limits_{\substack{p|d\\p\nmid N\\p \text{ prime}}} (1-p^{-s})\\
 &=&\frac{N^s d^s \prod\limits_{\substack{p|d \text{ or }p| N\\p \text{ prime}}} (1-p^{-s})}{N^s\prod\limits_{\substack{p| N\\p \text{ prime}}} (1-p^{-s})}\\
 &=&\frac{J_s(Nd)}{J_s(N)} = \frac{J_s(\frac{n}{t})}{J_s(\frac{n}{dt})} 
\end{eqnarray*}

\end{proof}

\end{thm}
We can now have some observations about the uniqueness of solutions if they exist. Note that in the above theorem, if  $d=1$, then obviously the number of solutions is 1. We see in the next corollary the other cases where uniqueness can be claimed.
\begin{cor}
 The restricted linear congruence (\ref{eq:singlevar}) has a unique solution only in the following two cases:
 \begin{enumerate}
  \item $d=1$
  \item $s=1, n=2^ru, t=2^{r-1}$ with $u,v$ odd and $v|u$ and $r\in\N$.
 \end{enumerate}
 \begin{proof}
 If $d=1$, then the number of solutions is clearly 1.
  Assume that there exists a prime $p$ such that $p|d, p\nmid N$. If no such $p$ exists, then the number of solutions is $d^s$ which is not equal to 1. Assume further that $d=p^j$ for some $j\geq 1$. By definition, $d^s=(a^{s},\frac{n^s}{t^s})_s$ (which is if and only if $d=(a,\frac{n}{t})$).  Let $p^r||n, p^k||t$. Then $0\leq j+k\leq r$. Now $p\nmid N\Rightarrow p\nmid p^{r-j-k} \Rightarrow r-j-k=0 \Rightarrow r=j+k$.
  So
  \begin{eqnarray*}
d^s\prod\limits_{\substack{p|d\\p\nmid N\\p \text{ prime}}} (1-p^{-s}) = 1 &\Rightarrow& p^{js}(1-p^{-s})=1\\
&\Rightarrow& p^{s(j-1)}(p^s-1)=1\\
&\Rightarrow& p^{s(j-1)}(p^s-1)=1\\
&\Rightarrow& p^{s(j-1)}=1 \text{ and }p^s-1=1\\
&\Rightarrow& s(j-1)=0 \text{ and }p^s=2\\
&\Rightarrow& s=1,p=2 \text{ and }s(j-1)=0\\
&\Rightarrow& s=1,p=2 \text{ and }j=1\\
  \end{eqnarray*}
Now the second statement also follows.

Note that from the above arguments, the only prime that can appear in $d$ is 2 and so the assumption that $d$ has only one prime factor covers the general case as well.
 \end{proof}
 \end{cor}
\section{Number of solutions of the linear congruence when $a_i=1$}
In this section we count the number of solutions of the linear congruence $x_1+\ldots+ x_k\equiv b(\text{mod }n^s)$ with restrictions $(x_i,n^s)_s=t_i^s$ where $t_i|n$ are given. A formula in this case was obtained in \cite{namboothiri4}, but for our purpose, we would like to derive a different formula. The following result generalizes \cite[Theorem 12, Theorem 12']{cohen1956extensionof} and \cite[Theorem 5.7]{namboothiri5}.
\begin{thm}
 Let $b\in\Z,n,t_i\in\N,t_i|n$ for $1\leq i\leq k$ given. The number of solutions of the linear congruence $x_1+\ldots+ x_k\equiv b(\text{mod }n^s)$ with restrictions $(x_i,n^s)_s=t_i^s$ is 
 \begin{equation}
N_{n,s}(b)=  \frac{1}{n^s}\sum\limits_{d|n}c_{d,s}(b)\prod\limits_{i=1}^{k}c_{\frac{n}{t_i},s}(\frac{n}{d})
 \end{equation}
\begin{proof}
 We use the DFT of $\rseven$ functions to prove our claim. Define 
 \begin{equation}
  \rho_{n,s,t}(m)=\begin{cases}
                   1 \text{ if }(m,n^s)_s=t^s\\
                   0 \text{ otherwise }
                  \end{cases}
 \end{equation}
Then $\rho_{n,s,t}$ is $\rseven$. So by the DFT of such functions, we have $\hat{\rho}_{n,s,t}(m)=\sum\limits_{d|n}\rho_{n,s,t}(d^s)c_{\frac{n}{d},s}(m)$. Note that $\rho_{n,s,t}(d^s)=1$  if and only if $(d^s,n^s)_s=t^s$. But since $d|n$, we have $(d^s,n^s)_s=d^s$. So $\rho_{n,s,t}(d^s)=1$  if and only if $d^s=t^s$ or $d=t$. Hence we have $\hat{\rho}_{n,s,t}(m) = c_{\frac{n}{t},s}(m)$. Consider the Cauchy product of the functions $\rho_{n,s,t_i}$ for $i=1,\ldots,k$ applied to $b$.
\begin{eqnarray*}
 (\rho_{n,s,t_1}\otimes \ldots\otimes \rho_{n,s,t_k})(b) &=& \sum\limits_{x_1+\ldots +x_k=b\,(\text{mod }n^s)}\rho_{n,s,t}(x_i)\ldots \rho_{n,s,t_k}(x_k)\\
 &=&  \sum\limits_{\substack{x_1+\ldots +x_k=b\,(\text{mod }n^s)\\(x_i,n^s)_s=t_i^s}} 1\\
 &=& N_{n,s}(b)
 \end{eqnarray*}
 So we get 
 \begin{eqnarray*}
  \widehat{N_{n,s}}(b) &=& \widehat{(\rho_{n,s,t_1}\otimes \ldots\otimes \rho_{n,s,t_k})}(b)\\
  &=&\prod\limits_{i=1}^k\hat{\rho}_{n,s,t_i}(b)\\
  &=&\prod\limits_{i=1}^k c_{\frac{n}{t_i},s}(b)
 \end{eqnarray*}
Now $\widehat{\widehat{N_{n,s}}}=n^sN_{n,s}$. So
\begin{eqnarray*}
 n^sN_{n,s}(b) &=& \widehat{\left(\prod\limits_{i=1}^k c_{\frac{n}{t_i},s}\right)}(b)\\
 &=& \sum\limits_{d|n}\left(\prod\limits_{i=1}^k c_{\frac{n}{t_i},s}(d)\right)c_{\frac{n}{d},s}(b)\\
 &=& \sum\limits_{d|n}c_{\frac{n}{d},s}(b)\left(\prod\limits_{i=1}^k c_{\frac{n}{t_i},s}(d)\right)
\end{eqnarray*}
which proves our claim.

\end{proof}

\end{thm}

\section{The general case}
We prove our main theorem here with the help of the two special cases  we had in the previous sections. The precise statement of our theorem is the following:
\begin{thm}
 Let $a_i,b\in\Z,n,s,t_i\in\N$ where $t_i|n$ be given. The number of solutions of the linear congruence 
 \begin{equation}\label{eq:gen_lin_cong2}
  a_1^{s}x_1+\ldots+a_k^{s} x_k \equiv b\,(\text{mod } n^s)
 \end{equation}
with restrictions $(x_i,n^s)_s=t_i^s$ is
\begin{equation}
 \frac{1}{n^s}\prod\limits_{i=1}^k\frac{J_s(\frac{n}{t_i})}{J_s(\frac{n}{t_id_i})}\sum\limits_{d|n}c_{d,s}(b)\prod\limits_{i=1}^k c_{\frac{n}{t_i d_i},s}(\frac{n}{d}) \text{ where } d_i=(a_i,\frac{n}{t_i})
\end{equation}
\begin{proof}
 Consider (\ref{eq:gen_lin_cong2}) and the congruence equation 
 \begin{equation}\label{eq:singlevar2}
y_1+\ldots+y_k \equiv b\,(\text{mod } n^s)  
 \end{equation}
If $\langle x_1,\ldots,x_k\rangle\in\Z^k$ is a solution for the first congruence, then $x_i$ is a solution for the congruence equation $a_i^{s}x_i=y_i \,(\text{mod } n^s)$. Conversely if $a_i^{s}x_i=y_i \,(\text{mod } n^s)$ has a solution $x_i$ for $i=1,\ldots,k$, then it gives a solution for (\ref{eq:gen_lin_cong2}). Now $a_i^{s}x_i=y_i \,(\text{mod } n^s) \Rightarrow (a_i^{s}x_i,n^s)_s=(y_i,n^s)_s$. Note that if solution $x_i$ with $(x_i,n^s)_s=t_i^s$ exists for the  congruence $a_i^{s}x_i=y_i \,(\text{mod } n^s)$  then $(a_i^{s}x_i,n^s)_s=(y_i,n^s)_s =t_i^sd_i^s$. The number of solutions for the congruence (\ref{eq:singlevar2}) with the restrictions $(y_i,n^s)_s =t_i^sd_i^s$ is $\frac{1}{n^s}\sum\limits_{d|n}c_{d,s}(b)\prod\limits_{i=1}^{k}c_{\frac{n}{t_id_i},s}(\frac{n}{d})$. Corresponding to each such solution $y_i$, we have $\frac{J_s(\frac{n}{t})}{J_s(\frac{n}{dt})} $ number of $x_i$'s. Hence the total number of solutions for (\ref{eq:gen_lin_cong2}) would be what we claimed.
\end{proof}

\end{thm}

\begin{rem}
 The most general case for the congruences we consider would be the one with the power $s$ in $a_i^s$ and $n^s$ removed while keeping the restrictions $(x_i,n)_s=t_i^s$. Let us consider the equation $ax\equiv b \,(\text{mod }n)$ and try to use the same proof given in theorem \ref{th:singlevar}. In the proof there, we used the fact that $(a^{s} \frac{x_0}{t^s}, \frac{n^s}{t^s})_s = (\frac{b}{t^s},\frac{n^s}{t^s})_s = d^s$
    and $(\frac{x_0}{t^s}, \frac{n^s}{t^s})_s = 1$ implies that $(a^{s}, \frac{n^s}{t^s})_s = d^s = (\frac{b}{t^s},\frac{n^s}{t^s})_s $. As we said, if we replace $a^s$ by $a$ and $n^s$ by $n$, then we have the conditions  $(a \frac{x_0}{t^s}, \frac{n}{t^s})_s = (\frac{b}{t^s},\frac{n}{t^s})_s = d^s$ and $(\frac{x_0}{t^s}, \frac{n}{t^s})_s = 1$. These two conditions need not give $(a, \frac{n}{t^s})_s = d^s$ as in the proof of theorem \ref{th:singlevar}.
\end{rem}

\bibliography{nt3} 

\begin{thebibliography}{10}

\bibitem{bibak2016restricted}
Khodakhast Bibak, Bruce~M Kapron, and Venkatesh Srinivasan.
\newblock On a restricted linear congruence.
\newblock {\em International Journal of Number Theory}, 12(08):2167--2171,
  2016.

\bibitem{bibak2017restricted}
Khodakhast Bibak, Bruce~M Kapron, Venkatesh Srinivasan, Roberto Tauraso, and
  L{\'a}szl{\'o} T{\'o}th.
\newblock Restricted linear congruences.
\newblock {\em Journal of Number Theory}, 171:128--144, 2017.

\bibitem{bibak2015almost}
Khodakhast Bibak, Bruce~M Kapron, Venkatesh Srinivasan, and L{\'a}szl{\'o}
  T{\'o}th.
\newblock On an almost-universal hash function family with applications to
  authentication and secrecy codes.
\newblock {\em arXiv preprint arXiv:1507.02331}, 2015.

\bibitem{brauer1926losung}
A~Brauer.
\newblock Lösung der aufgabe 30.
\newblock {\em Jber. Deutsch. Math.--Verein}, (35):92--94, 1926.

\bibitem{cohen1949extension}
Eckford Cohen.
\newblock An extension of ramanujan's sum.
\newblock {\em Duke Math. J}, 16(85-90):2, 1949.

\bibitem{cohen1955class}
Eckford Cohen.
\newblock A class of arithmetical functions.
\newblock {\em Proceedings of the National Academy of Sciences},
  41(11):939--944, 1955.

\bibitem{cohen1956extensionof}
Eckford Cohen.
\newblock An extension of ramanujan's sum. iii. connections with totient
  functions.
\newblock {\em Duke Math. J}, 23:623--630, 1956.

\bibitem{dixon1960finite}
J~D Dixon.
\newblock A finite analogue of the goldbach problem.
\newblock {\em Canadian Mathematical Bulletin}, 3:121--126, 1960.

\bibitem{lehmer1913certain}
D~N Lehmer.
\newblock Certain theorems in the theory of quadratic residues.
\newblock {\em The American Mathematical Monthly}, 20(5):151--157, 1913.

\bibitem{liskovets2010multivariate}
Valery~A Liskovets.
\newblock A multivariate arithmetic function of combinatorial and topological
  significance.
\newblock {\em Integers}, 10(1):155--177, 2010.

\bibitem{mccarthy1960generation}
Paul~J McCarthy.
\newblock The generation of arithmetical identities.
\newblock {\em J. reine angew. Math}, 203:55--63, 1960.

\bibitem{montgomery2006multiplicative}
Hugh~L Montgomery and Robert~C Vaughan.
\newblock {\em Multiplicative number theory I: Classical theory}, volume~97.
\newblock Cambridge University Press, 2006.

\bibitem{namboothiri5}
K~Vishnu Namboothiri.
\newblock The discrete fourier transform of $(r,s)$-even functions.
\newblock {\em arXiv preprint arXiv:1708.04507 [math.NT]}, 2017.

\bibitem{namboothiri4}
K~Vishnu Namboothiri.
\newblock On solving a restricted linear congruence using generalized ramanujan
  sums.
\newblock {\em arXiv preprint arXiv:1708.04505 [math.NT]}, 2017.

\bibitem{nicol1954sterneck}
Charles~A Nicol and Harry~S Vandiver.
\newblock A von sterneck arithmetical function and restricted partitions with
  respect to a modulus.
\newblock {\em Proceedings of the National Academy of Sciences},
  40(9):825--835, 1954.

\bibitem{radekacher1925aufgabe}
H~Rademacher.
\newblock Aufgabe 30.
\newblock {\em Jber. Deutsch. Math.--Verein}, 158(34), 1925.

\bibitem{sander2013adding}
JW~Sander and T~Sander.
\newblock Adding generators in cyclic groups.
\newblock {\em Journal of Number Theory}, 133(2):705--718, 2013.

\bibitem{terras1999fourier}
Audrey Terras.
\newblock {\em Fourier analysis on finite groups and applications}, volume~43.
\newblock Cambridge University Press, 1999.

\end{thebibliography}
\bibliographystyle{plain}

\end{document}